\documentclass[11pt]{amsart}
\usepackage{amsmath, amssymb,color, amsthm}
\usepackage[utf8]{inputenc}
\usepackage{bigints}
\usepackage{hyperref, xcolor, tikz}
\usepackage{enumitem, subcaption}
\usepackage{standalone}
\usepackage{longtable}
\usepackage{makecell}
\usepackage{etoolbox}
\usepackage[thinlines]{easytable}
\usepackage{booktabs,longtable,array,ragged2e}
\newcolumntype{P}[1]{>{\RaggedRight\hspace{0pt}}p{#1}}

\makeatletter
\patchcmd{\@maketitle}
  {\ifx\@empty\@dedicatory}
  {\ifx\@empty\@date \else {\vskip3ex \centering\footnotesize\@date\par\vskip1ex}\fi
   \ifx\@empty\@dedicatory}
  {}{}
\patchcmd{\@adminfootnotes}
  {\ifx\@empty\@date\else \@footnotetext{\@setdate}\fi}
  {}{}{}
\makeatother

\usepackage{geometry}
\usepackage{placeins}
\usepackage[normalem]{ulem}

\definecolor{goetheblau}{cmyk}{1.00 0.20 00 0.40}
\definecolor{hellgrau}{cmyk}{0.04 0.04 0.05 0.02}
\definecolor{sandgrau}{cmyk}{0.12 0.09 0.13 0}
\definecolor{dunkelgrau}{cmyk}{0.25 0.25 0.30 0.75}
\definecolor{purple}{cmyk}{0.08 1.00 0.30 0.36}
\definecolor{emorot}{cmyk}{0.04 1.00 0.80 0.07}
\definecolor{senfgelb}{cmyk}{0.01 0.25 1.00 0.05}
\definecolor{gruen}{cmyk}{0.62 0.40 0.87 0.09}
\definecolor{magenta}{cmyk}{0.08 0.86 0.12 0.12}
\definecolor{orange}{cmyk}{0 0.70 1.00 0.04}
\definecolor{sonnengelb}{cmyk}{0 0.12 0.95 0}
\definecolor{hellesgruen}{cmyk}{0.40 0.17 0.81 0.07}
\definecolor{lichtblau}{cmyk}{0.80 00 0.06 0.04}

\hypersetup{%
  colorlinks = true,
  linkcolor  = emorot,
  citecolor = hellesgruen
}

\newcommand{\Z}{\mathbb{Z}}

\newcommand{\Q}{\mathbb{Q}}

\newcommand{\F}{\mathbb{F}}

\usepackage{commath}
\usepackage{esdiff}
\usepackage[]{multirow}

\DeclareMathOperator{\Mat}{\text{Mat}}
\DeclareMathOperator{\Frob}{\text{Frob}}

\newcommand\restr[2]{{
  \left.\kern-\nulldelimiterspace 
  #1 
  \vphantom{\big|} 
  \right|_{#2} 
  }}

\usepackage{graphicx}
\usepackage{mathrsfs}

\usepackage[textwidth=20mm]{todonotes}
\usepackage{autonum}
\usetikzlibrary{arrows.meta,bending}
\usetikzlibrary{patterns}
\usepackage{pict2e,color}
\usepackage{placeins}
\title[Paramodular forms from Calabi-Yau operators]{Paramodular forms from
Calabi-Yau operators}

\begin{document}
\author[]{Nutsa Gegelia}
    \address{Nutsa Gegelia, Institut f\"ur Mathematik, Johannes Gutenberg-Universit\"at Mainz, Germany
}
\email{gegelian@uni-mainz.de}
\author[]{Duco van Straten}
    \address{Duco van Straten, Institut f\"ur Mathematik, Johannes Gutenberg-Universit\"at Mainz, Germany
}
\email{straten@mathematik@uni-mainz.de}
\maketitle
\noindent
\allowdisplaybreaks
\begin{abstract}
In this note we report on the conjectural identification of paramodular forms from Calabi-Yau motives of Hodge type $(1,1,1,1)$ of moderately low conductor. We calculate Euler factors from Calabi-Yau operators from the AESZ list
\cite{almkvist507430tables}, \cite{CYDB} by the method described in \cite{candelas2021local},
seek a fit with the  tables provided by \cite{assaf2023database} and for consistency check the
approximate functional equation for the Euler product for primes $<1000$ numerically, using
\cite{dokchitser2004computing}, \cite{hadoop}.
\end{abstract}

\section{Introduction}

Consider a threefold $X$ defined over $\Z$ and smooth over $\Q$. By Poincar\'e duality, the third
$\ell$--adic cohomology $H^3X:=H^3_{et}(\overline{X},\Q_\ell)$ carries a non-degenerate alternating form and
thus defines a {\em symplectic $Gal(\overline{\Q}/\Q)$ representation}. 
The {\em $L$--function} for the Galois representation $H^3X$ is defined in the usual way as product of
its Euler factors: 
\begin{equation}
    L(H^3X,s):=\prod_p E_p(p^{-s})^{-1},
\end{equation}
where the {\em Euler factor at the prime $p$} is determined  via
\begin{equation}
    E_p(T) := \det\left(1-T\cdot\Frob_p \mid H^3X^{I_p}\right),
\end{equation}
$I_p$ is the inertia subgroup at $p$ and $\Frob_p$ is a Frobenius element. The {\em conductor $N$} of
the Galois representation is defined as in \cite{serre1969facteurs}; the primes $p$ with $p\mid N$ are called bad primes and the ones with $p\nmid N$ good primes,
for which the reduction of $X \bmod p$ remains smooth.

If $X$ is a Calabi-Yau manifold, we have Hodge numbers $(h^{30},h^{21},h^{12}, h^{30})=(1,\tau,\tau,1)$,
where $\tau$ is the dimension of the local deformation space of $X$. In case $\tau=0$, we are dealing with
a {\em rigid Calabi-Yau manifold} and $H^3X$ is a two-dimensional Galois representation.
It was proven in \cite{gouvea2011rigid} that $X$ is {\em modular}, in the sense that there exists a $\mathbb{Q}$-Hecke eigenform $f$ of weight $4$ for some congruence subgroup $\Gamma_0(N)$ such that
\[ L(H^3X,s)=L(f,s),\]
but as far as we know, it is unknown if all such modular forms have such a rigid Calabi-Yau incarnation.
In \cite{meyer2005modular} one finds many examples; in many cases there are birationally different rigid Calabi-Yaus
giving the same modular form $f$.

In the case $\tau=1$ we are dealing with Hodge structures of type $(1,1,1,1)$ and corresponding four-dimensional
symplectic Galois representations. For good primes $p$ the Euler factor  $E_p(T)$ has the special form
\begin{equation}\label{eulerqq}
E_p(T) = 1 + \alpha_p T + \beta_p p T^2 + \alpha_p p^3 T^3 + p^6 T^4.
\end{equation}
For the bad primes the polynomial $E_p$ will have lower degree, but in general it is harder to obtain its precise form. In the simplest cases the polynomial is of degree 3 and factors in a linear and a quadratic factor.

For general reasons explained in \cite{gross2016langlands}, $L$-functions coming from such four-dimensional
Galois representations $H^3X$ and with conductor $N$ can be expected to be $L$-functions of Siegel modular
($\mathbb{Q}$-Hecke eigen-) forms of weight three on the three-dimensional Siegel space $\mathbb{H}_2$, invariant under the paramodular group $K(N)$  of level $N$, (paramodular form for short). So one may ask
\begin{equation}\label{paramodularity}
L(H^3X,s) \stackrel{?}{=} L(F,s), \;\;\; F \in S_3(\mathbb{H}_2/K(N)).
\end{equation}

Again the obvious question arises: which paramodular forms admit a realisation as $H^3X$ of a Calabi-Yau threefold $X$? Special cases arise if $H^3X$ is reducible, and $F$ is a lift, but in this note we are mainly interested in the generic cases, where $F$ is a non-lift.

The quest to exhibit matches of $L$-functions of Calabi-Yau threefolds with those of Siegel modular forms was explicitly adressed in \cite{cohen2015computing}, but was hampered at the time by lack of information on paramodular Hecke-eigenforms. Due to the hard work of many, this now has changed and
lists of paramodular forms have appeared in \cite{rama2020computation} and \cite{assaf2023database}. In \cite{bruinier2021moduli} and
\cite{golyshev2023congruences}  the identification of some Calabi-Yau spaces conjecturally realising forms of low conductor were reported on. The story of the paramodular form of lowest level $61$ will be treated in detail in a forthcoming paper.

In this note we report on Euler factors obtained from Calabi-Yau motives obtained from
Calabi-Yau operators from  \cite{almkvist507430tables} that match the Euler factors with those of  paramodular forms appearing in the table in \cite{assaf2023database}. Furthermore, by proposing Euler factors for the bad primes, we check the functional equation of the corresponding completed $L$-function
numerically, using the Pari implementation \cite{hadoop} of Dokchitser's algorithm for computing $L$-functions numerically
\cite{dokchitser2004computing}.

\section{Calabi-Yau operators}
Although there exists a large number of examples of Calabi-Yau threefolds of different
topologies, the number of examples with Hodge number $h^{12}=1$ is rather restricted.
As $h^{12}=\tau$ has the interpretation the dimension of local moduli space, these examples always appear in the form of one-parameter families, usually as pencils, 
i.e. as fibres $X_t$ of families $f: \mathcal{X} \longrightarrow {\mathbb{P}}^1$ defined over $\Q$. 
If the equations of $X_t$ are simple enough, one can determine Euler factors of $H^3X_t$ by
direct point counting or the method of Gauss sums. The 14 families corresponding to the mirror
families of complete intersections in weighted projective spaces can be studied in this way, \cite{CdOV}.
Other examples can be described in terms of fibre products of families of elliptic curves
and  the Euler factors may be obtained again via an explicit point-counting. But this direct method
becomes problematic in more complicated geometries. Many interesting Calabi-Yau varieties are
defined by Laurent polynomials, where direct counting is possible in principle,
but needs additional consideration of compactification and resolution to identify the pure
weight 3 part in the cohomology.\\

The idea here is to use the Picard-Fuchs equation $\mathcal{P} \in \Q \langle t,\frac{d}{dt}\rangle$
of the family $f: \mathcal{X} \longrightarrow {\mathbb{P}}^1$
describing the variation of cohomology groups $H^3X_t$. In many (but not all) examples there is a point of maximal unipotent monodromy (MUM-point), in which case $\mathcal{P}$ is a so-called {\em Calabi-Yau operator} in the sense of \cite{almkvist507430tables}. (For an overview see
also \cite{van2017calabi}).
It turns out there is a surprising method to calculate Euler factors of $H^3X_t$ directly from
the Calabi-Yau operator alone, without having to know the equations for the variety $X_t$.
This method builds on ideas from Dwork, and was developed further by Candelas, de la Ossa and van Straten
\cite{candelas2021local}.

\section{$p$-adic model of Frobenius}

We now give a brief sketch of this method, which remains conjectural in most cases.
A Calabi-Yau operator of order 4 with a MUM-point at $0$ 
is written in the form 
\begin{equation}
    \mathcal{P}=\theta^4+tP_1(\theta)+t^2P_2(\theta)+\ldots+t^rP_r(\theta) \in \mathbb{Q}\langle t,\theta\rangle,\;\;\theta=t\frac{d}{dt},
\end{equation}
where the $P_i$-s are polynomials of degree four. The number $r$ is called the {\em degree} of the
operator, so that the operators of degree $1$ are hypergeometric.
The operator $\mathcal{P}$ can be written also as
\begin{equation}
    \mathcal{P}=\Delta(t) \frac{d^4}{dt^4}+\ldots \in \mathbb{Q}\langle t,\frac{d}{dt}\rangle,
  \end{equation}
where the polynomial $\Delta(t)$ defines the singularities of the operator.
There is a Frobenius basis of solutions around 0 of the form
\begin{equation}\label{frobeniusbasis}
\begin{aligned}
    f_0 &= A(t), \\
    f_1 &= A(t)\log(t) + B(t),\\
    f_2 &= \frac{1}{2}A(t)\log(t)^2 + B(t)\log(t) + C(t),\\
    f_3 &= \frac{1}{6}A(t)\log(t)^3 + \frac{1}{2}B(t)\log(t)^2+C(t)\log(t)+D(t), 
\end{aligned}
\end{equation}
where $A(t)\in\Q[|t|],B(t),C(t),D(t)\in t\Q[|t|]$ are power series.  Let
\begin{equation}
    E(t) := \begin{pmatrix}
        A & \theta(A) & \theta^2(A) & \theta^3(A) \\
        B & A + \theta(B) & 2\theta(A) + \theta^2(B) & 3\theta^2(A) + \theta^3(B)\\
        C & B + \theta(C) & A + 2\theta(B) + \theta^2(C) & 3\theta(A) + 3\theta^2(B) + \theta^3(C) \\
        D & C + \theta(D) & B + 2\theta(C) + \theta^2(D) & 3\theta(B) + 3\theta^2(C) + \theta^3(D) 
    \end{pmatrix}.
  \end{equation}
  Following Dwork \cite{dwork1969padiccycles}, we define for a prime $p$ the \emph{U-matrix}
  \begin{equation}\label{umatrixorder2}
    U_p(t) := E(t^p)^{-1}\cdot U_p(0)\cdot E(t) \in \Mat(4,4,\Q[[t]]).
\end{equation}
For the  \emph{limit U-matrix} $U_p(0)$ we make the following {\em Ansatz}:
\begin{equation}
    U_p(0):=\begin{pmatrix}
        1 & 0 & 0 & 0 \\
        0 & p & 0 & 0 \\
        0 & 0 & p^2 & 0 \\
        p^3\cdot x_p & 0 & 0 & p^3
    \end{pmatrix},
  \end{equation}
where 
\begin{equation}
x_p := x\cdot\zeta_p(3),
\end{equation}
and $x \in \Q$ only depends on the operator, but not on $p$. For the appropriate choice of $x$
the matrix series $U_p(t)$ is $p$-adic integral and has an expansion
\begin{equation}\label{umatrixisrational}
    U_p(t) = \frac{V_0(t)}{\Delta(t)^{p\cdot\delta_0}} + p\cdot \frac{V_1(t)}{\Delta(t)^{p\cdot\delta_1}} +  p^2\frac{V_2(t)}{\Delta(t)^{p\cdot\delta_2}} + p^3\cdot \frac{V_3(t)}{\Delta(t)^{p\cdot\delta_3}} + p^4\cdot \frac{V_4(t)}{\Delta(t)^{p\cdot\delta_4}} + p^5 \cdot \frac{V_5(t)}{\Delta(t)^{p\cdot\delta_5}} + \dots.
\end{equation}
Here the $V_i(t)\in \Mat(4,4,\Q[t])$ are polynomial matrices of predictable degree $d_i$ and $\Delta(t)$ as
defined above. The exponents $\delta_n$ depend on the operator in question. In many simple
cases (and the appropriate choice of $x$) one has $\delta_0=\delta_1=\delta_2=\delta_3=0$ or $\delta_0=\delta_1=\delta_2=0$, $\delta_3=1$.

Conjecturally one has, possibly up to twist by a character:
\begin{equation}
U_p(\mu(t)) \stackrel{?}{=} \text{Matrix of} \;  \Frob_p:H^3 X_t\rightarrow H^3 X_t
\end{equation}
where $\mu(t) \in \Z_p$ is the Teichm\"uller lift of $t \in \F_p$.
So each Euler factor \eqref{eulerqq} of $H^3X_t$ is equal to
\begin{equation}
E_{p,t}(T) =\det(1- U(\mu(t)) T).
\end{equation}
As the reciprocal roots of $E_p(T)$ are supposed to have size $p^{3/2}$ (Weil bounds), it is enough (for $p \ge 7$) to compute \eqref{umatrixisrational} up to $\bmod p^4$
and lift the coefficient $\bmod p^4$ to $\mathbb{Z}$, so that the zeros fulfill the Weil bounds.\\

\section{The Euler factory}
In the AESZ list \cite{almkvist507430tables} one finds a collection of Calabi-Yau operators, which by definition
posses a point of maximal unipotent monodromy. These
operators are all of geometrical origin, and define a family of 'four-dimensional Calabi-Yau motives' which in good cases 
come from nice pencils of Calabi-Yau threefolds. In \cite{CYDB} one finds an online database presenting the operators from \cite{almkvist507430tables}, complemented
with operators found later, ordered by their degree $r$.
The symbol $a.b$ stands for the $b$th operator of degree $a$
in the database,  so that $1.1$ is the hypergeometric Picard-Fuchs operator from the mirror quintic pencil.

Using the above sketched algorithm and this list of operators, one in fact has an {\em Euler factory:}
picking an operator $\mathcal{P}$, a (regular) value for $t$ and a (good) prime $p$, one can compute an Euler factor.
The computation performed with this method is very fast and can be performed for a large range of primes $p$.
Up to now most calculations were performed on a 64 bit, Intel 12th Gen Core i5 4.4GHz, 16Gb RAM without specific parallelization. They can be carried on up to $p<1000$ within a few hours and could be done without much effort for higher $p$. The calculations
could easily be parallelized and performed on a cluster, which would improve performance significantly.
The main restricting factor is the lack of larger memory.

The database \cite{assaf2023database} contains information on Hecke eigenvalues $a_p$ and $b_p$ for $7587$ paramodular Hecke-eigenforms of generic type (G). Of these there are 2575 forms defined over $\mathbb{Q}$. After running a search by matching the coefficients $\alpha_p$ with the eigenvalues $a_p$ of paramodular forms in this list, we were able to find more than 300 matches. These are  listed in Table \ref{matches1}. 
One may observe that different cases may lead to the same paramodular form. An example 
is the case of the operators 5.22 and 5.25 at $t=1$ corresponding to the same paramodular form 2.K.73.3.0.a.a.
In such a case we say that the paramodular form has different incarnations. Standard conjectures of Hodge-Tate type predict that there exists a  correspondence relating the two Calabi-Yau geometries.

In total we have 233 singletons, 32 doublets and 3 triplets, covering a total of 268 paramodular forms.

\subsection*{Checking the functional equation}
As a sanity check of the $U$-matrix method, one would also like to make a numerical check of the functional equation
of the product of Euler factors. We did this for some cases (Tables \ref{listindatabase1} and \ref{listnotindatabase}). For this we used the PARI/GP package for $L$-functions \cite{hadoop}, which implements
the algorithms of Dokchitser \cite{dokchitser2004computing}. 

For Hodge numbers $(1,1,1,1)$, the completed $L$--function $\Lambda(s)$ is defined as
\begin{equation}\label{functeq}
\Lambda(s)= \left(\frac{N}{\pi^4}\right)^{s/2}\Gamma\left(\frac{s-1}{2}\right)\Gamma\left(\frac{s}{2}\right)\Gamma\left(\frac{s}{2}\right)\Gamma\left(\frac{s+1}{2}\right)L(H^3X,s),
\end{equation}
where $N$ is the conductor of the Galois representation. The
functional equation then reads:
\begin{equation}\label{functionaleq}
    \Lambda(s)=\epsilon\Lambda(4-s),
\end{equation}
where $\epsilon=\pm$ is the {\em sign} of the functional equation. To check this, we need sufficiently many Euler factors, including the Euler factors for all the primes of bad reduction. Remarkably, it turns out that with some care, the $U$-matrix method produces in many cases also the correct Euler factors at primes of bad reduction. However, there
remain many cases where we have to guess the 'right' Euler factors at the bad primes. In the Appendix there are tables of Euler factors given as triples $[p,\alpha_p,\beta_p]$ for each case from the Tables \ref{listindatabase1} and \ref{listnotindatabase} together with the corresponding bad Euler factors.
With Philip Candelas and Xenia de la Ossa we plan to convert our heuristics into more precise formulae for the
conductor $N$ and the Euler factors for the bad primes \cite{CGOvS}.

\begin{figure}
\begin{center}
\includegraphics[height=8cm]{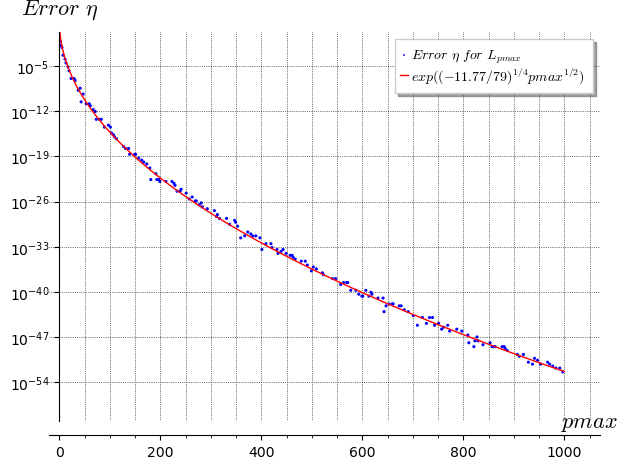}
\caption{Precision curve for 79a.}
\label{fig:figure1}
\end{center}
\end{figure}

With this data available, we consider the partial Euler product 
\begin{equation}
    L_{pmax}(H^3X,s)=\prod_{p=2}^{pmax} E_p(p^{-s})^{-1}
\end{equation}
and we use the function \verb|lfuncheckfeq| from PARI/GP  \cite{hadoop} for $\Lambda(s)$. This function returns a bit accuracy $b$ such that the numerical error for $L(s)$ is smaller than $2^b$. For more details see the User's guide of PARI/GP. We determine the {\em precision},
which we define as the reciprocal of the numerical error $\eta$ in \verb|lfuncheckfeq|.
We consider $\eta$ as a function of the prime number $pmax$. It turns out that $\eta$ as function of $pmax$ has a
universal form described by the formula
\[\eta(pmax)= \exp\left(-\frac{c}{N^{1/4}}\cdot pmax^{1/2}\right),\;\;\; c \approx  11.77.\]
By tracking these precision curves one can conveniently monitor the increase of precision: if not every Euler factor up to
$pmax$ is correct, the precision stagnates after the wrong factor.
We illustrate this typical behaviour for the form 79a in Figure \ref{fig:figure1}.

We consider this as good evidence that the $\alpha_p$ up to 1000 and the $\beta_p$ up to at least 31 are correct.
Note that the error $\eta$ is not monotonically falling, but rather exhibits typical accidental drops to higher
precision than the formula would predict.

\begin{center}
\begin{figure}[t]
\includegraphics[height=8cm]{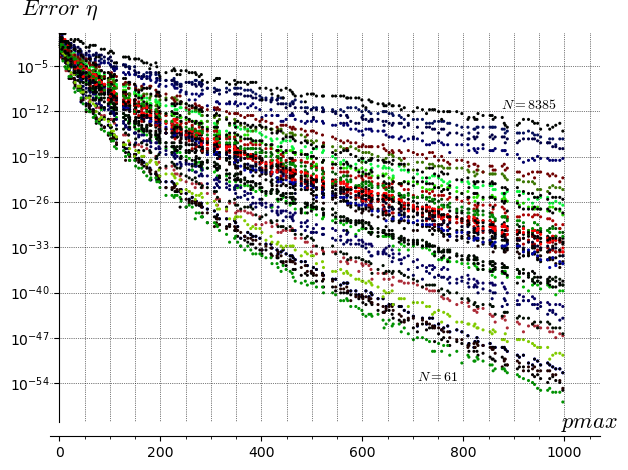}
\caption{A sample of precision curves drawn in one plot; It can be observed that the precision grows approximately in the following way:
$
    \eta(N,pmax) \sim \exp\left(-\frac{c}{N^{1/4}}\cdot pmax^{1/2}\right), \; c\approx 11.77.
$
}
\label{fig:figure2}
\end{figure}
\end{center}

Figure \ref{fig:figure2} shows the striking universality of the precision curves we obtained from the examples appearing in the Tables \ref{listindatabase1} and \ref{listnotindatabase}.

\subsection*{Congruences}
Another good reason to compute many complete Euler factors, rather than only Frobenius traces which would suffice for
the purpose of identification, is that it allows us to get better statistics on factorisation properties. If for a
certain prime $\ell$ and all good $p$ the polynomial $E_p(T)\; \mod\; \ell$ is reducible, we say we have an
$\ell$-congruence. The primes and the corresponding decomposition types are indicated in the rightmost columns  of the Tables \ref{listindatabase1} and \ref{listnotindatabase}.
Of course, it may happen that for specific $t$ all Euler factors factor over $\mathbb{Z}$. Such points are so-called {\em attractor points} and are also of great interest, but these cases are not included here, as the corresponding paramodular form is lift, and not of generic type (G).

\subsection*{Faltings-Serre method}
A conjectural match of $L$-functions
\[L(H^3X,s)=L(F,s)\]
can be proven by establishing equality of sufficiently many Euler factors. For two-dimensional Galois
representations the Faltings-Serre method provides a computable bound that can be used to make this strategy effective. In \cite{brumer2019paramodularity} this was extended to four-dimensional $GSp_4$ representations, which has found ample application. However, so far method presupposes that the corresponding $\bmod 2$ Galois
representation is irreducible. As almost all matches in our list carry a 2-congruence, this method is currently
applicable only in a few cases: 2.K.61.3.0.a.a, 2.K.94.3.0.a.b, 2.K.302.3.0.a.a, 2.K.431.3.0.a.b. Details on these examples will appear elsewhere.

\subsection*{Predictions for higher conductors}
In most cases the Euler factory will produce an $L$-function with conductor
$N >1000$, which brings us outside the range of the Table \ref{matches1}. Nevertheless, based on some
heuristics for the conductor $N$, we were able to verify the functional equation for examples of rather high
conductor. Table \ref{listnotindatabase} contains a sample, far from complete, of such examples with conductor
ranging between $1000$ and $10.000$, but we were able to do examples with conductor as high as $\sim 200.000$.
For all these cases we predict the existence of a paramodular form of the indicated conductor. So in this sense one may
consider the Euler factory as an {\em oracle}, that predicts the existence of certain paramodular forms.

We also found an example with low conductor 128 not appearing 
in \cite{assaf2023database}, but which appears in \cite{rama2020computation}:
\small
\begin{table}[h]
\hspace*{-3em}
\centering
\renewcommand{\arraystretch}{1.3}
\begin{tabular}[h]{c|c|c|c|c|c|c|c}
Corresponds to form & Conductor & \cite{CYDB} & \cite{almkvist507430tables} & $t$ & Precision & $\epsilon$ in \eqref{functionaleq} & Type of congruence mod $\ell$ \\
\cline{2-8}
in \cite{rama2020computation} & 128 & $2.55$ & \#42 & 1/16 & E-47 & + & 2, $(1^4)$\\
\end{tabular}
\end{table}
\normalsize

\subsection*{Further calculations}
We reported on more than $300$ potential identifications of paramodular forms from Calabi-Yau operators, which represent
about $10 \%$ of the tabulated forms. But we have no doubts by refining our search we will be able to identify many more.
As we already remarked, the $U$-matrix method is easily parallelisable and we plan to start a comprehensive large scale
implementation of the Euler factory and seek correspondences between Calabi-Yau motives and paramodular forms. For this
we seek to cooperate with Albrecht Klemm and Paul Blesse, who also have identified some of the above reported matches.

\subsection*{Acknowledgments}
We would like to thank Philip Candelas, Neil Dummigan, Vasily Golyshev, Abhiram Kidambi, Xenia de la Ossa, Ariel Pacetti, Gustavo Rama,
Gonzalo Tornar\'ia, John Voight, Wadim Zudilin and other participants of the working group {\em GdT} for
showing interest in the subject and useful feedback on earlier drafts of this text.
Special thanks to the Abdus Salam ICTP for hospitality and scientific exchange during the Workshop/School on Number Theory and Physics in June 2024.  We especially thank Jan Lukas Igelbrink for generous help with programming issues.
We gratefully acknowledge the support by the Deutsche Forschungsgemeinschaft (DFG, German Research Foundation) – Project-ID 444845124 – TRR 326.

\bibliography{quellen} 

\begin{thebibliography}{CYDB}
\expandafter\ifx\csname url\endcsname\relax
  \def\url#1{\texttt{#1}}\fi
\expandafter\ifx\csname doi\endcsname\relax
  \def\doi#1{\burlalt{doi:#1}{http://dx.doi.org/#1}}\fi
\expandafter\ifx\csname urlprefix\endcsname\relax\def\urlprefix{URL }\fi
\expandafter\ifx\csname href\endcsname\relax
  \def\href#1#2{#2}\fi
\expandafter\ifx\csname burlalt\endcsname\relax
  \def\burlalt#1#2{\href{#2}{#1}}\fi

\bibitem[CYDB]{CYDB}
P.~Metelitsyn and D.~van Straten.
\newblock Calabi-{Y}au database.
\newblock \urlprefix\url{https://cydb.mathematik.uni-mainz.de/}.

\end{thebibliography}


\newcommand{\etalchar}[1]{$^{#1}$}
\begin{thebibliography}{CdlORV00}
\expandafter\ifx\csname url\endcsname\relax
  \def\url#1{\texttt{#1}}\fi
\expandafter\ifx\csname doi\endcsname\relax
  \def\doi#1{\burlalt{doi:#1}{http://dx.doi.org/#1}}\fi
\expandafter\ifx\csname urlprefix\endcsname\relax\def\urlprefix{URL }\fi
\expandafter\ifx\csname href\endcsname\relax
  \def\href#1#2{#2}\fi
\expandafter\ifx\csname burlalt\endcsname\relax
  \def\burlalt#1#2{\href{#2}{#1}}\fi

\bibitem[AESZ]{almkvist507430tables}
G.~Almkvist, C.~van Enckevort, D.~van Straten, and W.~Zudilin.
\newblock Tables of {C}alabi--{Y}au equations, 2010,
  \burlalt{math/0507430}{http://arxiv.org/abs/math/0507430}.
\newblock \urlprefix\url{https://arxiv.org/abs/math/0507430}.

\bibitem[ALR{\etalchar{+}}24]{assaf2023database}
E.~Assaf, W.~Ladd, G.~Rama, G.~Tornar\'ia, and J.~Voight.
\newblock A database of paramodular forms from quinary orthogonal modular
  forms.
\newblock In {\em Lu{C}a{NT}: {LMFDB}, computation, and number theory}, volume
  796 of {\em Contemp. Math.}, pages 243--259. Amer. Math. Soc., [Providence],
  RI, [2024] \copyright2024.

\bibitem[BPP{\etalchar{+}}19]{brumer2019paramodularity}
A.~Brumer, A.~Pacetti, C.~Poor, G.~Tornar\'ia, J.~Voight, and D.~S. Yuen.
\newblock On the paramodularity of typical abelian surfaces.
\newblock {\em Algebra Number Theory}, 13(5):1145--1195, 2019.

\bibitem[CdlORV00]{CdOV}
P.~Candelas, X.~de~la Ossa, and F.~Rodriguez-Villegas.
\newblock Calabi-{Y}au manifolds over finite fields, i, 2000,
  \burlalt{hep-th/0012233}{http://arxiv.org/abs/hep-th/0012233}.
\newblock \urlprefix\url{https://arxiv.org/abs/hep-th/0012233}.

\bibitem[CdlOvS21]{candelas2021local}
P.~Candelas, X.~de~la Ossa, and D.~van Straten.
\newblock Local zeta functions from {C}alabi-{Y}au differential equations,
  2021, \burlalt{2104.07816}{http://arxiv.org/abs/2104.07816}.
\newblock \urlprefix\url{https://arxiv.org/abs/2104.07816}.

\bibitem[CGdlOvS]{CGOvS}
P.~Candelas, N.~Gegelia, X.~de~la Ossa, and D.~van Straten.
\newblock Work in progress.

\bibitem[Coh15]{cohen2015computing}
H.~Cohen.
\newblock Computing {$L$}-functions: a survey.
\newblock {\em J. Th\'eor. Nombres Bordeaux}, 27(3):699--726, 2015.

\bibitem[CYDB]{CYDB}
P.~Metelitsyn and D.~van Straten.
\newblock Calabi-{Y}au database.
\newblock \urlprefix\url{https://cydb.mathematik.uni-mainz.de/}.

\bibitem[Dok04]{dokchitser2004computing}
T.~Dokchitser.
\newblock Computing special values of motivic {$L$}-functions.
\newblock {\em Experiment. Math.}, 13(2):137--149, 2004.

\bibitem[Dwo69]{dwork1969padiccycles}
B.~Dwork.
\newblock {$p$}-adic cycles.
\newblock {\em Inst. Hautes \'Etudes Sci. Publ. Math.}, (37):27--115, 1969.

\bibitem[Gro16]{gross2016langlands}
B.~K. Gross.
\newblock On the {L}anglands correspondence for symplectic motives.
\newblock {\em Izv. Ross. Akad. Nauk Ser. Mat.}, 80(4):49--64, 2016.

\bibitem[GvS23]{golyshev2023congruences}
V.~Golyshev and D.~van Straten.
\newblock Congruences via fibered motives.
\newblock {\em Pure Appl. Math. Q.}, 19(1):233--265, 2023.

\bibitem[GY11]{gouvea2011rigid}
F.~Q. Gouv\^ea and N.~Yui.
\newblock Rigid {C}alabi-{Y}au threefolds over {$\Bbb Q$} are modular.
\newblock {\em Expo. Math.}, 29(1):142--149, 2011.

\bibitem[Mey05]{meyer2005modular}
C.~Meyer.
\newblock {\em Modular {C}alabi-{Y}au threefolds}, volume~22 of {\em Fields
  Institute Monographs}.
\newblock American Mathematical Society, Providence, RI, 2005.

\bibitem[{PAR}]{hadoop}
{PARI/GP L-functions package}.
\newblock \urlprefix\url{https://people.maths.bris.ac.uk/~matyd/computel/}.

\bibitem[RT20]{rama2020computation}
G.~Rama and G.~Tornar\'ia.
\newblock Computation of paramodular forms.
\newblock In {\em A{NTS} {XIV}---{P}roceedings of the {F}ourteenth
  {A}lgorithmic {N}umber {T}heory {S}ymposium}, volume~4 of {\em Open Book
  Ser.}, pages 353--370. Math. Sci. Publ., Berkeley, CA, 2020.

\bibitem[Ser69]{serre1969facteurs}
J.-P. Serre.
\newblock Facteurs locaux des fonctions z{\^e}ta des vari{\'e}t{\'e}s
  alg{\'e}briques (d{\'e}finitions et conjectures).
\newblock {\em S{\'e}minaire Delange-Pisot-Poitou. Th{\'e}orie des nombres},
  11(2):1--15, 1969.

\bibitem[vS18]{van2017calabi}
D.~van Straten.
\newblock Calabi-{Y}au operators.
\newblock In {\em Uniformization, {R}iemann-{H}ilbert correspondence,
  {C}alabi-{Y}au manifolds \& {P}icard-{F}uchs equations}, volume~42 of {\em
  Adv. Lect. Math. (ALM)}, pages 401--451. Int. Press, Somerville, MA, 2018.

\bibitem[vS21]{bruinier2021moduli}
D.~van Straten.
\newblock Rank four {C}alabi-{Y}au motives of low conductor.
\newblock In {\em Oberwolfach report: Moduli spaces and Modular forms (hybrid
  meeting)}. Bruinier, Jan Hendrik and van der Geer, Gerard and Gritsenko,
  Valery, 2021.

\end{thebibliography}
\bibliographystyle{habbrv}

\small
\begin{table}[!ht]
\caption{All the matches found so far}\label{matches1}
  \centering
 \begin{minipage}{0.44\linewidth} 
    \begin{tabular}{|c|c|c|c|}
    \hline
\cite{assaf2023database} & \cite{CYDB} & $t$ & Cond. \\
    \hline
2.K.61.3.0.a.a  & 5.24  & 1     & 61  \\
2.K.69.3.0.a.a  & 5.39  & 1/32  & 69  \\
2.K.73.3.0.a.a  & 5.22  & 1     & 73  \\
2.K.73.3.0.a.a  & 5.25  & 1     & 73  \\
2.K.76.3.0.a.a  & 5.27  & --1    & 76  \\
2.K.76.3.0.a.a  & 5.28  & --1/32 & 76  \\
2.K.79.3.0.a.a  & 2.5   & --1    & 79  \\
2.K.87.3.0.a.a  & 5.15  & --1/64 & 87  \\
2.K.87.3.0.a.a  & 5.31  & 1     & 87  \\
2.K.89.3.0.a.a  & 5.18  & 1     & 89  \\
2.K.94.3.0.a.b  & 4.34  & --1    & 94  \\
2.K.96.3.0.a.a  & 5.6   & 1/8   & 96  \\
2.K.104.3.0.a.a & 6.22  & 4     & 104 \\
2.K.104.3.0.a.a & 8.40  & --1/3  & 104 \\
2.K.112.3.0.a.a & 8.41  & --1/8  & 112 \\
2.K.112.3.0.a.a & 8.42  & --1/16 & 112 \\
2.K.116.3.0.a.a & 5.7   & --1    & 116 \\
2.K.116.3.0.a.a & 5.46  & 1     & 116 \\
2.K.118.3.0.a.a & 2.69  & --1/32 & 118 \\
2.K.119.3.0.a.a & 2.53  & 1     & 119 \\
2.K.123.3.0.a.a & 2.54  & 1/81  & 123 \\
2.K.129.3.0.a.a & 8.22  & 1     & 129 \\
2.K.129.3.0.a.a & 8.23  & 1     & 129 \\
2.K.130.3.0.a.b & 2.62  & --1    & 130 \\
2.K.138.3.0.a.a & 5.110 & 1/16  & 138 \\
2.K.140.3.0.a.a & 8.2   & --1/8  & 140 \\
2.K.153.3.0.a.a & 2.53  & --1    & 153 \\
2.K.153.3.0.a.a & 2.54  & --1/81 & 153 \\
2.K.160.3.0.a.a & 2.52  & --1/64 & 160 \\
2.K.161.3.0.a.a & 2.55  & 1     & 161 \\
2.K.162.3.0.a.a & 8.66  & 1/64  & 162 \\
2.K.162.3.0.a.a & 8.67  & 1     & 162 \\
2.K.165.3.0.a.a & 6.9   & 1     & 165 \\
2.K.167.3.0.a.a & 8.16  & --1    & 167 \\
2.K.172.3.0.a.b & 5.7   & 1     & 172 \\
2.K.172.3.0.a.b & 5.46  & --1    & 172 \\
2.K.172.3.0.a.a & 5.84  & 1     & 172 \\
2.K.173.3.0.a.a & 8.16  & 1     & 173 \\
2.K.176.3.0.a.a & 2.69  & 1/16  & 176 \\
2.K.182.3.0.a.a & 1.4   & 1     & 182 \\
2.K.184.3.0.a.b & 5.18  & 1/2   & 184 \\
2.K.197.3.0.a.a & 5.24  & --1    & 197 \\
2.K.203.3.0.a.a & 5.30  & 1     & 203 \\
2.K.205.3.0.a.a & 1.6   & --1    & 205 \\
2.K.208.3.0.a.c & 3.15  & 1/8   & 208 \\
2.K.209.3.0.a.b & 8.6   & --1    & 209 \\
2.K.213.3.0.a.a & 8.20  & --1    & 213 \\
2.K.224.3.0.a.b & 2.52  & 1/8   & 224 \\
2.K.224.3.0.a.a & 2.55  & 1/4   & 224 \\
2.K.232.3.0.a.c & 2.57  & 1/64  & 232 \\
2.K.238.3.0.a.b & 2.64  & --1    & 238 \\
    \hline
    \end{tabular}
 \end{minipage} 
 \begin{minipage}{0.36\linewidth}
    \begin{tabular}{|c|c|c|c|}
    \hline
\cite{assaf2023database} & \cite{CYDB} & $t$ & Cond. \\
    \hline 
2.K.239.3.0.a.a & 10.1  & --1    & 239 \\
2.K.239.3.0.a.a & 11.3  & --1/4  & 239 \\
2.K.245.3.0.a.a & 5.39  & --1/32 & 245 \\
2.K.248.3.0.a.c & 1.5   & --1/64 & 248 \\
2.K.255.3.0.a.b & 1.3   & 1     & 255 \\
2.K.257.3.0.a.a & 1.3   & --1    & 257 \\
2.K.262.3.0.a.a & 8.28  & --1    & 262 \\
2.K.262.3.0.a.a & 8.29  & --1/16 & 262 \\
2.K.264.3.0.a.b & 5.49  & --1/4  & 264 \\
2.K.264.3.0.a.e & 6.17  & 1/3   & 264 \\
2.K.269.3.0.a.a & 5.30  & --1    & 269 \\
2.K.275.3.0.a.a & 8.1   & --1    & 275 \\
2.K.275.3.0.a.a & 8.1   & --1/8  & 275 \\
2.K.281.3.0.a.a & 10.1  & 1     & 281 \\
2.K.281.3.0.a.a & 11.3  & --1/2  & 281 \\
2.K.284.3.0.a.a & 8.59  & --1    & 284 \\
2.K.287.3.0.a.a & 4.39  & --1    & 287 \\
2.K.295.3.0.a.a & 5.15  & 1/64  & 295 \\
2.K.295.3.0.a.a & 5.31  & --1    & 295 \\
2.K.296.3.0.a.a & 2.69  & 1/64  & 296 \\
2.K.296.3.0.a.a & 5.95  & 1/16  & 296 \\
2.K.296.3.0.a.b & 8.22  & 1/2   & 296 \\
2.K.296.3.0.a.b & 8.23  & 2     & 296 \\
2.K.297.3.0.a.g & 2.56  & 1/81  & 297 \\
2.K.302.3.0.a.a & 4.34  & 1     & 302 \\
2.K.304.3.0.a.a & 5.79  & 1/2   & 304 \\
2.K.312.3.0.a.a & 5.27  & --2    & 312 \\
2.K.312.3.0.a.a & 5.28  & --1/64 & 312 \\
2.K.315.3.0.a.b & 2.52  & 1     & 315 \\
2.K.316.3.0.a.a & 7.12  & --1/4  & 316 \\
2.K.320.3.0.a.d & 2.65  & 1/16  & 320 \\
2.K.320.3.0.a.a & 7.1   & --1/8  & 320 \\
2.K.327.3.0.a.b & 2.61  & --1    & 327 \\
2.K.327.3.0.a.c & 11.18 & 1     & 327 \\
2.K.328.3.0.a.a & 5.23  & --1    & 328 \\
2.K.328.3.0.a.b & 5.76  & --1/2  & 328 \\
2.K.330.3.0.a.a & 7.3   & 1/8   & 330 \\
2.K.335.3.0.a.a & 8.86  & 1     & 335 \\
2.K.335.3.0.a.a & 11.12 & --1/5  & 335 \\
2.K.336.3.0.a.a & 8.69  & 1/4   & 336 \\
2.K.337.3.0.a.a & 9.6   & --1    & 337 \\
2.K.341.3.0.a.a & 5.5   & 1     & 341 \\
2.K.341.3.0.a.a & 5.16  & 1/32  & 341 \\
2.K.344.3.0.a.a & 5.70  & --1/4  & 344 \\
2.K.348.3.0.a.a & 6.26  & --1/2  & 348 \\
2.K.353.3.0.a.a & 2.55  & --1    & 353 \\
2.K.357.3.0.a.a & 2.60  & 1     & 357 \\
2.K.359.3.0.a.a & 2.56  & --1    & 359 \\
2.K.363.3.0.a.b & 8.65  & 1     & 363 \\
2.K.364.3.0.a.b & 5.1   & 1     & 364 \\
2.K.368.3.0.a.h & 1.5   & 1/64  & 368 \\
\hline
    \end{tabular}
    \end{minipage}
\end{table}
\normalsize

\small
\begin{table}[!ht]
\ContinuedFloat
\caption{All the matches found so far (continued)}
  \centering
 \begin{minipage}{0.44\linewidth} 
    \begin{tabular}{|c|c|c|c|}
    \hline
\cite{assaf2023database} & \cite{CYDB} & $t$ & Cond. \\
    \hline
2.K.381.3.0.a.c & 11.16 & --1    & 381 \\
2.K.384.3.0.a.e & 1.3   & 1/64  & 384 \\
2.K.384.3.0.a.b & 3.10  & 1/4   & 384 \\
2.K.390.3.0.a.h & 2.64  & 1     & 390 \\
2.K.391.3.0.a.a & 11.18 & --1    & 391 \\
2.K.396.3.0.a.b & 3.1   & --1/2  & 396 \\
2.K.419.3.0.a.a & 8.33  & --1/8  & 419 \\
2.K.419.3.0.a.a & 8.34  & 1     & 419 \\
2.K.424.3.0.a.b & 2.57  & 1/16  & 424 \\
2.K.431.3.0.a.c & 1.5   & 1     & 431 \\
2.K.431.3.0.a.a & 2.5   & 1     & 431 \\
2.K.431.3.0.a.b & 2.6   & --1    & 431 \\
2.K.433.3.0.a.a & 1.5   & --1    & 433 \\
2.K.435.3.0.a.a & 5.49  & --1    & 435 \\
2.K.436.3.0.a.a & 8.63  & 1/8   & 436 \\
2.K.442.3.0.a.a & 5.81  & 1     & 442 \\
2.K.448.3.0.a.f & 1.5   & --1/16 & 448 \\
2.K.448.3.0.a.o & 1.11  & --1/64 & 448 \\
2.K.448.3.0.a.c & 2.55  & 1/64  & 448 \\
2.K.448.3.0.a.b & 8.43  & 1/8   & 448 \\
2.K.455.3.0.a.f & 2.65  & 1     & 455 \\
2.K.456.3.0.a.b & 1.3   & --2    & 456 \\
2.K.456.3.0.a.d & 5.102 & 1/4   & 456 \\
2.K.462.3.0.a.b & 2.69  & --1    & 462 \\
2.K.462.3.0.a.c & 6.21  & --1    & 462 \\
2.K.465.3.0.a.a & 5.4   & 1     & 465 \\
2.K.469.3.0.a.a & 8.35  & --1/2  & 469 \\
2.K.475.3.0.a.a & 8.9   & --1    & 475 \\
2.K.476.3.0.a.f & 3.1   & 1/8   & 476 \\
2.K.476.3.0.a.d & 8.5   & 1/8   & 476 \\
2.K.480.3.0.a.b & 1.3   & 1/16  & 480 \\
2.K.480.3.0.a.i & 3.34  & --1/16 & 480 \\
2.K.482.3.0.a.a & 8.85  & 1/64  & 482 \\
2.K.483.3.0.a.b & 2.53  & 1/9   & 483 \\
2.K.486.3.0.a.b & 3.4   & --1/27 & 486 \\
2.K.488.3.0.a.c & 2.56  & --1/16 & 488 \\
2.K.491.3.0.a.b & 5.18  & --1    & 491 \\
2.K.492.3.0.a.a & 5.75  & --1/8  & 492 \\
2.K.492.3.0.a.c & 5.85  & 1     & 492 \\
2.K.495.3.0.a.a & 5.4   & --1    & 495 \\
2.K.496.3.0.a.c & 2.62  & 1/2   & 496 \\
2.K.498.3.0.a.c & 5.82  & --1    & 498 \\
2.K.502.3.0.a.a & 8.24  & 1/8   & 502 \\
2.K.503.3.0.a.d & 2.56  & 1     & 503 \\
2.K.503.3.0.a.b & 11.16 & 1     & 503 \\
2.K.504.3.0.a.a & 2.53  & 8     & 504 \\
2.K.504.3.0.a.e & 4.39  & --2    & 504 \\
2.K.504.3.0.a.n & 5.124 & 8/9   & 504 \\
2.K.506.3.0.a.a & 5.70  & 1     & 506 \\
2.K.510.3.0.a.e & 2.62  & 1/81  & 510 \\
2.K.519.3.0.a.a & 5.68  & 1     & 519 \\
    \hline
    \end{tabular}
 \end{minipage}
 \begin{minipage}{0.36\linewidth}
    \begin{tabular}{|c|c|c|c|}
    \hline
\cite{assaf2023database} & \cite{CYDB} & $t$ & Cond. \\
    \hline 
2.K.520.3.0.a.f & 3.1   & --1    & 520 \\
2.K.521.3.0.a.a & 5.51  & --1    & 521 \\
2.K.532.3.0.a.d & 2.2   & --1/8  & 532 \\
2.K.532.3.0.a.a & 8.47  & --1/8  & 532 \\
2.K.535.3.0.a.a & 2.61  & 1     & 535 \\
2.K.540.3.0.a.d & 5.111 & 3/4   & 540 \\
2.K.542.3.0.a.b & 5.79  & --1    & 542 \\
2.K.544.3.0.a.h & 1.3   & --1/16 & 544 \\
2.K.544.3.0.a.c & 2.1   & --1/8  & 544 \\
2.K.544.3.0.a.k & 2.55  & 1/2   & 544 \\
2.K.544.3.0.a.o & 5.32  & 1/8   & 544 \\
2.K.552.3.0.a.a & 5.42  & 1/64  & 552 \\
2.K.554.3.0.a.c & 5.20  & 1     & 554 \\
2.K.558.3.0.a.g & 3.9   & --1/3  & 558 \\
2.K.560.3.0.a.h & 2.64  & 1/2   & 560 \\
2.K.560.3.0.a.g & 6.12  & 1/2   & 560 \\
2.K.561.3.0.a.d & 2.65  & --1/81 & 561 \\
2.K.572.3.0.a.a & 2.10  & 1     & 572 \\
2.K.583.3.0.a.a & 5.2   & 1     & 583 \\
2.K.585.3.0.a.b & 2.52  & 1/25  & 585 \\
2.K.585.3.0.a.a & 8.35  & 1/2   & 585 \\
2.K.586.3.0.a.b & 5.70  & --1    & 586 \\
2.K.588.3.0.a.j & 6.21  & --1/3  & 588 \\
2.K.588.3.0.a.j & 7.10  & 1/9   & 588 \\
2.K.599.3.0.a.b & 5.35  & 1     & 599 \\
2.K.603.3.0.a.h & 2.20  & 1/81  & 603 \\
2.K.603.3.0.a.a & 8.63  & 1     & 603 \\
2.K.605.3.0.a.c & 5.2   & --1    & 605 \\
2.K.610.3.0.a.a & 8.47  & --1/32 & 610 \\
2.K.612.3.0.a.h & 5.125 & --1/8  & 612 \\
2.K.615.3.0.a.d & 2.54  & --1/25 & 615 \\
2.K.621.3.0.a.d & 2.55  & 1/3   & 621 \\
2.K.621.3.0.a.a & 8.43  & 1/3   & 621 \\
2.K.623.3.0.a.b & 3.25  & 1     & 623 \\
2.K.624.3.0.a.j & 5.79  & --1/2  & 624 \\
2.K.624.3.0.a.h & 6.26  & 1     & 624 \\
2.K.630.3.0.a.e & 5.81  & --1    & 630 \\
2.K.632.3.0.a.a & 5.23  & 1/2   & 632 \\
2.K.634.3.0.a.b & 8.24  & --1/8  & 634 \\
2.K.640.3.0.a.j & 1.3   & --1/64 & 640 \\
2.K.640.3.0.a.s & 2.18  & 1/32  & 640 \\
2.K.640.3.0.a.d & 2.52  & --1/16 & 640 \\
2.K.640.3.0.a.h & 2.64  & 1/32  & 640 \\
2.K.640.3.0.a.u & 3.5   & --1/4  & 640 \\
2.K.640.3.0.a.t & 5.123 & --1/64 & 640 \\
2.K.640.3.0.a.u & 6.38  & 1/8   & 640 \\
2.K.644.3.0.a.f & 3.1   & 1/2   & 644 \\
2.K.649.3.0.a.a & 4.36  & --1    & 649 \\
2.K.650.3.0.a.b & 5.77  & 1     & 650 \\
2.K.656.3.0.a.f & 2.54  & 1/16  & 656 \\
2.K.656.3.0.a.b & 8.73  & --1/10 & 656 \\
\hline
    \end{tabular}
    \end{minipage}
\end{table}
\normalsize

\small
\begin{table}[!ht]
\ContinuedFloat
\caption{All the matches found so far (continued)}
  \centering
 \begin{minipage}{0.44\linewidth} 
    \begin{tabular}{|c|c|c|c|}
    \hline
\cite{assaf2023database} & \cite{CYDB} & $t$ & Cond. \\
    \hline
2.K.663.3.0.a.c & 2.57  & --1/9  & 663 \\
2.K.672.3.0.a.h & 1.3   & 1/4   & 672 \\
2.K.672.3.0.a.o & 1.6   & 1/16  & 672 \\
2.K.672.3.0.a.d & 5.11  & --1/16 & 672 \\
2.K.672.3.0.a.o & 8.55  & 1/2   & 672 \\
2.K.672.3.0.a.o & 8.55  & 1/8   & 672 \\
2.K.679.3.0.a.a & 9.6   & 1     & 679 \\
2.K.680.3.0.a.b & 2.62  & --1/16 & 680 \\
2.K.688.3.0.a.c & 5.95  & --1/16 & 688 \\
2.K.688.3.0.a.e & 11.5  & 1/16  & 688 \\
2.K.693.3.0.a.d & 5.102 & 1/27  & 693 \\
2.K.696.3.0.a.i & 5.82  & 2     & 696 \\
2.K.696.3.0.a.f & 8.26  & 1     & 696 \\
2.K.696.3.0.a.f & 8.27  & --1/32 & 696 \\
2.K.700.3.0.a.b & 8.73  & --1/7  & 700 \\
2.K.703.3.0.a.b & 8.4   & --1/8  & 703 \\
2.K.704.3.0.a.t & 2.56  & 1/4   & 704 \\
2.K.704.3.0.a.b & 5.73  & 1/16  & 704 \\
2.K.720.3.0.a.v & 2.62  & 1/4   & 720 \\
2.K.720.3.0.a.l & 9.5   & --1/4  & 720 \\
2.K.720.3.0.a.l & 11.9  & --1/8  & 720 \\
2.K.725.3.0.a.a & 8.7   & --1    & 725 \\
2.K.730.3.0.a.d & 1.4   & --1    & 730 \\
2.K.731.3.0.a.a & 6.9   & --1    & 731 \\
2.K.735.3.0.a.i & 3.33  & --1    & 735 \\
2.K.735.3.0.a.h & 4.71  & 1     & 735 \\
2.K.740.3.0.a.e & 2.11  & --1/16 & 740 \\
2.K.741.3.0.a.a & 6.26  & 1/8   & 741 \\
2.K.741.3.0.a.a & 8.51  & 1/6   & 741 \\
2.K.745.3.0.a.b & 5.51  & 1     & 745 \\
2.K.745.3.0.a.a & 8.20  & 1     & 745 \\
2.K.750.3.0.a.f & 2.64  & 1/25  & 750 \\
2.K.750.3.0.a.d & 3.9   & 1     & 750 \\
2.K.764.3.0.a.a & 9.6   & --1/2  & 764 \\
2.K.764.3.0.a.a & 11.6  & --1/3  & 764 \\
2.K.765.3.0.a.a & 5.55  & 1/9   & 765 \\
2.K.768.3.0.a.a & 5.120 & 1/64  & 768 \\
2.K.776.3.0.a.b & 8.24  & 1/4   & 776 \\
2.K.780.3.0.a.e & 5.43  & 1/16  & 780 \\
2.K.801.3.0.a.a & 4.39  & 1     & 801 \\
2.K.806.3.0.a.b & 2.69  & 1     & 806 \\
2.K.806.3.0.a.e & 5.27  & 1/32  & 806 \\
2.K.806.3.0.a.e & 5.28  & 1     & 806 \\
2.K.806.3.0.a.a & 5.87  & --1    & 806 \\
2.K.806.3.0.a.d & 11.20 & --1    & 806 \\
2.K.808.3.0.a.b & 8.16  & --1/2  & 808 \\
2.K.812.3.0.a.c & 5.92  & --1/4  & 812 \\
2.K.814.3.0.a.d & 8.37  & 1     & 814 \\
2.K.814.3.0.a.d & 8.38  & 1/81  & 814 \\
2.K.816.3.0.a.b & 5.75  & 1/4   & 816 \\
2.K.816.3.0.a.d & 9.10  & --1/2  & 816 \\
    \hline
    \end{tabular}
 \end{minipage} 
 \begin{minipage}{0.36\linewidth}
    \begin{tabular}{|c|c|c|c|}
    \hline
\cite{assaf2023database} & \cite{CYDB} & $t$ & Cond. \\
    \hline 
2.K.832.3.0.a.c & 5.95  & --1/8  & 832 \\
2.K.832.3.0.a.k & 8.46  & 1/16  & 832 \\
2.K.832.3.0.a.b & 8.65  & 1/2   & 832 \\
2.K.840.3.0.a.h & 6.24  & 3/4   & 840 \\
2.K.840.3.0.a.h & 7.6   & 1/6   & 840 \\
2.K.840.3.0.a.h & 7.8   & 1/10  & 840 \\
2.K.840.3.0.a.g & 9.2   & 1/4   & 840 \\
2.K.847.3.0.a.c & 4.37  & --1/9  & 847 \\
2.K.848.3.0.a.a & 8.86  & 1/2   & 848 \\
2.K.848.3.0.a.a & 11.12 & --1/4  & 848 \\
2.K.856.3.0.a.a & 5.84  & --1/2  & 856 \\
2.K.860.3.0.a.e & 5.85  & --1    & 860 \\
2.K.867.3.0.a.f & 2.63  & 1/81  & 867 \\
2.K.867.3.0.a.g & 6.20  & 1/3   & 867 \\
2.K.867.3.0.a.g & 8.54  & --1/5  & 867 \\
2.K.867.3.0.a.g & 8.82  & 1/5   & 867 \\
2.K.870.3.0.a.a & 5.104 & 1     & 870 \\
2.K.877.3.0.a.a & 8.22  & --1    & 877 \\
2.K.877.3.0.a.a & 8.23  & --1    & 877 \\
2.K.880.3.0.a.g & 5.49  & 1/4   & 880 \\
2.K.882.3.0.a.h & 3.9   & --1/7  & 882 \\
2.K.884.3.0.a.b & 7.10  & 1     & 884 \\
2.K.891.3.0.a.h & 2.56  & --1/9  & 891 \\
2.K.896.3.0.a.g & 2.24  & 1/16  & 896 \\
2.K.896.3.0.a.o & 2.55  & 1/8   & 896 \\
2.K.896.3.0.a.m & 2.55  & 1/32  & 896 \\
2.K.899.3.0.a.a & 5.68  & --1    & 899 \\
2.K.908.3.0.a.a & 8.58  & 1     & 908 \\
2.K.910.3.0.a.b & 3.8   & --1/81 & 910 \\
2.K.912.3.0.a.i & 2.69  & 1/8   & 912 \\
2.K.928.3.0.a.a & 5.35  & 1/4   & 928 \\
2.K.928.3.0.a.a & 8.65  & --1/4  & 928 \\
2.K.930.3.0.a.a & 5.110 & --1/16 & 930 \\
2.K.945.3.0.a.k & 1.3   & 1/81  & 945 \\
2.K.945.3.0.a.t & 2.9   & 1/9   & 945 \\
2.K.945.3.0.a.a & 5.47  & --1/9  & 945 \\
2.K.951.3.0.a.a & 5.78  & --1/9  & 951 \\
2.K.952.3.0.a.p & 1.4   & 8     & 952 \\
2.K.952.3.0.a.o & 2.55  & 2     & 952 \\
2.K.952.3.0.a.n & 3.1   & 2     & 952 \\
2.K.960.3.0.a.h & 2.9   & 1/64  & 960 \\
2.K.960.3.0.a.e & 2.61  & 1/8   & 960 \\
2.K.960.3.0.a.x & 2.13  & --1/16 & 960 \\
2.K.962.3.0.a.c & 5.86  & --1    & 962 \\
2.K.975.3.0.a.b & 2.60  & --1    & 975 \\
2.K.976.3.0.a.c & 11.14 & --1    & 976 \\
2.K.981.3.0.a.a & 5.111 & --1    & 981 \\
2.K.987.3.0.a.b & 2.65  & --1    & 987 \\
2.K.992.3.0.a.l & 3.26  & --1/2  & 992 \\
2.K.992.3.0.a.e & 5.88  & 1/16  & 992 \\
2.K.996.3.0.a.d & 5.87  & --1/8  & 996 \\
\hline
    \end{tabular}
    \end{minipage}
\end{table}
\normalsize


\small
\begin{table}[h]
\caption{Matches with checked functional equation\label{listindatabase1}}
\centering
\renewcommand{\arraystretch}{1.3}
\begin{tabular}[h]{|c|c|c|c|c|c|c|c|}
\hline
Conductor &\cite{CYDB} & \cite{almkvist507430tables} & $t$ & \cite{assaf2023database} &  Precision & $\epsilon$ in \eqref{functionaleq} & Type of congruence mod $\ell$ \\
\hline
61 & $5.24$ & \#195& 1 & 2.K.61.3.0.a.a & E-57  & + & 19, (1,1,2); 43, (1,1,2) \\
69 & $5.39$ & \#224 & 1/32 & 2.K.69.3.0.a.a & E-57 & + & 2, (1,1,2); 19, (1,1,2)\\
73 & $5.22$ & \#193& 1 & 2.K.73.3.0.a.a & E-55 & + & 2, (2,2); 3, $(1^2, 1^2)$; 13, (1,1,2) \\
79 & $2.5$ & $\#25/A*b$ & --1 & 2.K.79.3.0.a.a & E-53 & + & 2, (1,1,2); 5, (2,2) \\
94 & $4.34$ & \#99 & --1 & 2.K.94.3.0.a.b & E-50 & + & 3, (1,1,2)\\
119 & $2.53$ & \#29 & 1 & 2.K.119.3.0.a.a & E-47 & + & 2, $(1^4)$; 3, $(1^2,1^2)$ \\
130 & $2.62$ & \#28 & --1 & 2.K.130.3.0.a.b & E-48 & + & 2, $(1^4)$; 3, (1,1,2); 7, (1,1,2)  \\
153 & $2.53$ & \#29 & --1 & 2.K.153.3.0.a.a & E-41 & + & 2, $(1^4)$; 3, $(1^2,1^2)$ \\
161 & $2.55$ & \#42 & 1 & 2.K.161.3.0.a.a & E-44 & + & 2, $(1^4)$ \\
182 & $1.4$ & $\#4/B*B$ & --1 & 2.K.182.3.0.a.a & E-43 & + & 2, $(1^2,2)$; 3, $(1^2,1^2)$ \\
197 & $5.24$ & \#195 & --1 & 2.K.197.3.0.a.a & E-41 & -- & 2, (1,1,2) \\
205 & 1.6 & \#6 & --1 & 2.K.295.3.0.a.a & E-41 & + & 2, $(1^4)$\\
224 & $2.55$ & \#42 & 1/4 & 2.K.224.3.0.a.a & E-40 & + & 2, $(1^4)$ \\
224 & $2.52$ & \#16 & 1/8 & 2.K.224.3.0.a.b & E-41 &  + & 2, $(1^4)$; 3, (1,1,2) \\
238 & $2.64$ & \#182 & --1 & 2.K.238.3.0.a.b & E-39 & + & 2, $(1^4)$; 3, (1,1,2); 11, (1,1,2) \\
245 & $5.39$ & \#224 & --1/32 & 2.K.245.3.0.a.a & E-38 & -- &  2, (1,1,2); 3, (2,2)\\
248 & 1.5 & $\#5/A*B$ & --1/64 & 2.K.248.3.0.a.c & E-40 & + & 2, (1,1,2); 3, (1,1,2)\\
255 & $1.3$ & $\#3/A*A$ & 1 & 2.K.255.3.0.a.b & E-39 & + & 2, $(1^4)$ \\
257 & $1.3$ & $\#3/A*A$ & --1 & 2.K.257.3.0.a.a & E-39 & + & 2, $(1^4)$ \\
302 & $4.34$ & \#99 & 1 & 2.K.315.3.0.a.b & E-36 & -- & None \\
315 & $2.52$ & \#16 & 1 & 2.K.302.3.0.a.a & E-37 & + & 2, $(1^4)$; 3, (1,1,2)  \\
327 & $2.61$ & \#26 & --1 & 2.K.327.3.0.a.b & E-36 & -- & 2, $(1^2,2)$ \\
353 & $2.55$ & \#42 & --1 & 2.K.353.3.0.a.a & E-36 & + & 2, $(1^4)$ \\
357 & $2.60$ & \#18 & 1 & 2.K.357.3.0.a.a & E-35 & -- & 2, $(1^4)$ \\
359 & $2.56$ & \#185 & --1 & 2.K.359.3.0.a.a & E-36 & + & 2, $(1^2,2)$; 3, (1,1,2)  \\
368 & 1.5 & $\#5/A*B$ & 1/64 & 2.K.368.3.0.a.h & E-35 & + & 2, (1,1,2); 3, (1,1,2)\\
384 & 1.3 & $\#3/A*A$ & 1/64 & 2.K.384.3.0.a.e & E-35 & + & 2, $(1^4)$\\
390 & $2.64$ & \#182 & 1 & 2.K.390.3.0.a.h & E-35 & + & 2, $(1,1,2)$; 3, (1,1,2); 11, (1,1,2) \\
396 & $3.1$ & \#34 & --1/2 & 2.K.396.3.0.a.b & E-36 & + & 2, $(1^4)$; 3, (1,1,2); 5, (1,1,2) \\
431 & $2.5$ & $\#25/A*b$ & 1 & 2.K.431.3.0.a.a & E-34 & -- & 2, (1,1,2); 5, (2,2)  \\
431 & $2.6$ & $\#24/B*b$ & --1 & 2.K.431.3.0.a.b & E-33 & -- & 3, (1,1,2); 5, (2,2)  \\
431 & $1.5$ & $\#5/A*B$ & 1 & 2.K.431.3.0.a.c & E-34 & + & 2, $(1^2,2)$; 3 $(1^2,2)$ \\
433 & $1.5$ & $\#5/A*B$ &--1 & 2.K.433.3.0.a.a & E-34 & + & 2, $(1^2,2)$; 3 $(1^2,2)$ \\
448 & $1.5$ & $\#5/A*B$ & --1/16 & 2.K.448.3.0.a.f & E-34 & + & 2, $(1^4)$; 3, (1,1,2)\\
448 & $1.11$ & \#11 & --1/64 & 2.K.448.3.0.a.o & E-34 & + &  2, $(1^4)$; 3, $(1^2,1^2)$\\
456 & $1.3$ & $\#3/A*A$ & --2 & 2.K.456.3.0.a.b & E-34 & + &  2, $(1^4)$; 3, (1,1,2)\\
462 & $2.69$ & \#205 & --1 & 2.K.462.3.0.a.b & E-33 & + & 2, $(1^4)$; 5, (1,1,2) \\
476 & $3.1$ & \#34 & 1/8 & 2.K.476.3.0.a.f & E-34 & + & 2, $(1^4)$; 3, (1,1,2); 5, (1,1,2) \\
480 & $1.3$ & $\#3/A*A$ & 1/16 & 2.K.480.3.0.a.b & E-34 & + &  2, $(1^4)$\\
496 &  $2.62$ & \#28 & 1/2 & 2.K.496.3.0.a.c & E-34 & + & 2, $(1^4)$; 3, (1,1,2); 7, (1,1,2)  \\
\hline
\end{tabular} 
\end{table}
\normalsize 

\small
\begin{table}[h]
\ContinuedFloat
\caption{Matches with checked functional equation (continued)}
\centering
\renewcommand{\arraystretch}{1.3}
\begin{tabular}[h]{|c|c|c|c|c|c|c|c|c|}
\hline
Conductor & \cite{CYDB} & \cite{almkvist507430tables} & $t$ &\cite{assaf2023database} &  Precision & $\epsilon$ in \eqref{functionaleq} & Type of congruence mod $\ell$ \\
\hline
503 & $2.56$ & \#185 & 1 & 2.K.503.3.0.a.d & E-33 & + & 2, $(1^2,2)$; 3, (1,1,2)  \\
520 & $3.1$ & \#34 & --1 & 2.K.520.3.0.a.f & E-34 & + & 2, $(1^4)$; 3, (1,1,2); 5, (1,1,2) \\
535 & $2.61$ & \#26 & 1 & 2.K.535.3.0.a.a & E-31 & -- & 2, $(1^2,2)$ \\
544 & $1.3$ & $\#3/A*A$ & --1/16 & 2.K.544.3.0.a.h & E-32 & + & 2, $(1^4)$ \\
544 & $2.55$ & \#42 & 1/2 & 2.K.544.3.0.a.k & E-32 & + & 2, $(1^4)$ \\
572 & $2.10$ & $\#70/B*c$ & 1 & 2.K.572.3.0.a.a & E-33 & + & 2, $(1^2,2)$; 3, $(1^2,1^2)$ \\
621 & $2.55$ & \#42 & 1/3 & 2.K.621.3.0.a.d & E-32 & + & 2, $(1^4)$ \\
640 & $1.3$ & $\#3/A*A$ & --1/64 & 2.K.640.3.0.a.j & E-31 & + & 2, $(1^4)$ \\
640& 2.52 & \#16 & $-1/16$ & 2.K.640.3.0.a.d & E-31 & + &  2, $(1^4)$, 3, (1,1,2)\\
644 & $3.1$ & \#34 & 1/2 & 2.K.644.3.0.a.f & E-31 & + & 2, $(1^4)$; 3, (1,1,2); 5, (1,1,2) \\
672 & $1.3$ & $\#3/A*A$ & 1/4 & 2.K.672.3.0.a.h & E-31 & + & 2, $(1^4)$ \\
672 & $1.6$ & \#6 & 1/16 & 2.K.672.3.0.a.o & E-32 & + & 2, $(1^4)$\\
704 & $2.56$ & \#185 & 1/4 & 2.K.704.3.0.a.t & E-29 & + & 2, $(1^2,2)$; 3, (1,1,2)  \\
730 & $1.4$ & $\#4/B*B$ & 1 & 2.K.730.3.0.a.d & E-30 & + & 2, $(2,2)$; 3, $(1^2,1^2)$ \\
806 & $2.69$ & \#205 & 1 & 2.K.806.3.0.a.b & E-28 & -- & 2, $(1^4)$; 5, (1,1,2) \\
891 & $2.56$ & \#185 & --1/9 & 2.K.891.3.0.a.h & E-28 & + & 2, $(1^2,2)$; 3, (1,1,2) \\
896 & $2.55$ & \#42 & 1/8 & 2.K.896.3.0.a.o & E-28 & + & 2, $(1^4)$ \\
945 & $1.3$ & $\#3/A*A$ & 1/81 & 2.K.945.3.0.a.k & E-28 & + & 2, $(1^4)$\\
952 & $1.4$ & $\#4/B*B$ & 8 & 2.K.952.3.0.a.p & E-28 & + & 2, (2,2); 3, $(1^2,1^2)$\\
952 & $3.1$ & \#34 & 2 & 2.K.952.3.0.a.n & E-28 & + & 2, $(1^4)$, 3, (1,1,2); 5, (1,1,2) \\
975 & $2.60$ & \#18 & --1 & 2.K.975.3.0.a.b & E-26 & -- & 2, $(1^4)$
\\
\hline
\end{tabular} 
\end{table}
\normalsize

\small
\begin{table}[h]
\caption{Candidates with checked functional equation (not in database \cite{assaf2023database})\label{listnotindatabase}}
\centering
\renewcommand{\arraystretch}{1.3}
\begin{tabular}[h]{|c|c|c|c|c|c|c|c|}
\hline
Conductor & \cite{CYDB} & \cite{almkvist507430tables} & $t$ & Precision & Sign in \eqref{functionaleq} & Type of congruence mod $\ell$ \\
\hline
1025 & $2.6$ & $\#24/B*b$ & 1 & E-26 & -- & 3, (1,1,2); 5, (2,2)  \\
1105 & 2.52 & $\#16$ & --1 & E-27 & + & 2, $(1^4)$; 3, (1,1,2)  \\
1562 & $1.1$ & $\#1/Quintic$ & 1 & E-23 & -- & 5, (1,1,1,1) \\
1584 & 2.62 & $\#28$ & --1/2 & E-25 & + & 2, $(1^4)$; 3, (1,1,2); 7, (1,1,2)  \\
1708 & $2.10$ & $\#70/B*c$ & --1 & E-25 & + & 2, $(1^2,2)$; 3, $(1^2,1^2)$ \\
1935 & $2.1$ & $\#45/A*a$ & --1 & E-22 & -- & 2, $(1^4)$; 3, (2,2)  \\
2145 & $2.9$ & $\#58/A*c$ & 1 & E-23 & + & 2, $(1^4)$; 3, (2,2)  \\
2159 & $2.1$ & $\#45/A*a$ & 1 & E-22 & -- & 2, $(1^4)$; 3, (2,2) \\
2465 & $2.9$ & $\#58/A*c$  & --1 & E-22 & + & 2, $(1^4)$; 3, (2,2); 5, (1,1,2) \\
2667 & $2.13$ & $\#36/A*d$& 1 & E-21 & + & 2, $(1^4)$ \\
3010 & $2.2$ & $\#15/B*a$ & 1 & E-20 & -- & 2, (2,2); 3, $(1^2,1^2)$  \\
3126 & $1.1$ & $\#1/Quintic$ & --1 & E-20 & + & 5, (1,1,1,1) \\
3391 & $2.7$ & $\#51/C*b$ & --1 & E-19 & -- & 2, $(1^2,2)$; 5, (2,2)  \\
4799 & $2.7$ & $\#51/C*b$ & 1 & E-20 & + & 2, $(1^2,2)$; 5, (2,2) \\
5642 & $2.2$ & $\#15/B*a$ & --1 & E-16 & -- & 2, (2,2); 3, $(1^2,1^2)$  \\
7057 & $2.20$ & $\#133/A*f$& --1 & E-17 & + & 2, $(1^2,2)$; 3, (2,2)  \\
8385 & $2.13$ & $\#36/A*d$ & --1 & E-16 & -- & 2, $(1^4)$ \\
\hline
\end{tabular} 
\end{table}
\normalsize

\end{document}